\documentclass[12pt]{article} 
\usepackage{amsfonts}
\usepackage{amssymb}
\usepackage{enumerate}
\usepackage{epsfig}
\usepackage{amsmath}
\date{\small\today}
\title{
 A Topology-Preserving Level Set Method for\\Shape Optimization}
\author{Oleg Alexandrov and Fadil Santosa\\
  \small University of Minnesota\\
  \small School of Mathematics \\
 \small \{aoleg, santosa\}@math.umn.edu}

\begin{document}
\maketitle
\begin{abstract}
The classical level set method, which represents the boundary of the
unknown geometry as the zero-level set of a function, has been shown
to be very effective in solving shape optimization problems.  The
present work addresses the issue of using a level  set representation
when there are simple geometrical and topological constraints.  We
propose a logarithmic barrier penalty which acts to enforce the
constraints, leading to an approximate solution to shape design
problems.
\end{abstract}

\section{Introduction}
The level set method
\cite{osher-sethian-1988,sethian-book,osher-fedkew-book} is a very
powerful approach for problems involving geometry and geometric
evolution.  It has also been applied to solving shape optimization
problems \cite{allaire, sethian-wiegmann, osher-santosa}, and it is at
this type of problems that this work is aimed.

By a \emph{shape} we mean a bounded region $D$ in $\mathbb R^n$
with $C^1$ boundary. One associates with $D$ a function $\phi:\mathbb
R^n\to\mathbb R$ with the property that $D$ is the \emph{level set} of
$\phi,$
$$
  D=\{x:\phi(x)>0\}.
$$
One then manipulates $D$ implicitly, through its level set
function $\phi$.  It is typical in shape optimization problems to start
with an initial shape, which is then improved in an iterative process.
Thus, one would start with a level set function $\phi(x)$ which is
updated at each iteration.

The advantage of the level set method is that it is much easier to
work with a globally defined function than to keep track of the
boundary of a domain. The latter, which can be achieved by using
marker points and spline interpolation, can become especially complicated if
$D$ has either several connected components, or is otherwise connected but has
several holes.  During the optimization process, the 
components or holes may merge or split, or even entirely disappear.
The level set method, on the other hand, takes care of
this kind of changes with great ease.

Given the shape $D$ there exist of course many functions $\phi$ whose
level set is $D$.  The most convenient $\phi$ to work with is the
\emph{signed distance} to the boundary $\partial D$ of $D,$ thus
\begin{equation}\label{1}
\phi(x)=
\begin{cases} 
  \mbox{dist}(x, \partial D), & x\in D,\\
  -\mbox{dist}(x, \partial D), & x\not \in D.
\end{cases}  
\end{equation}
Then $\phi$ will have the additional property 
\begin{equation}\label{2}
  \nabla\phi(x)\cdot\nabla\phi(x)=1
\end{equation}
for $x$ in a neighborhood of $\partial D.$ Any level set function
$\phi$ can be reinitialized as the signed distance to the set
$\{x:\phi(x)=0\},$ so from here on we will assume that $\phi$ always
satisfies \eqref{1}, by reinitializing it if necessary.

It is very easy to describe deformations of $D$ in terms of its level set
function $\phi$. For example, if $h:\mathbb R^n\to\mathbb R$ is a
function with $\sup|h(x)|$ small enough, then the level set of
$\phi+h$ is obtained from the level set $D$ of $\phi$ by shifting
every point $x\in \partial D$ by approximately the amount $h(x)$ in
the direction of the external normal to $\partial D$ at $x$ (which is
$-\nabla \phi(x)$).

While the level set method has its strong points -- one being that it
gives a representation that is topology-independent -- it is not obvious
how to extend it to problems where there are constraints.  Simple volume
(area in 2-D) constraints are relatively easy to incorporate
\cite{osher-santosa}.  Other constraints, such as a bound on the size of
a connected component of $D$, or the requirement that $D$ has a fixed number of
connected components, are not as easy to handle.  It is towards this class of
problems that this work is directed.

Our approach starts with the concept of subdomain neighborhood.  The
neighborhood of one subdomain will detect the nearness of other
subdomains, and will thus allow us to take action to prevent geometry or
topology changes.  This strategy can be formulated as a penalty
functional, which we describe in the next section.  We illustrate this
method by two numerical examples in Section 3.

We wish to mention the paper \cite{han} which also suggests a way of
adapting the level set method to preserve topology. The authors of this
paper do it in the context of image segmentation.  The key difference
between our work and \cite{han} is that their method is pixel-based.
The algorithm in \cite{han} is able to detect that a shape is about to
change topology only when certain dimensions of the shape are of size
comparable to the grid size. In the context of image processing this makes
a lot of sense, as then it is convenient to define a body to be
connected as long as it is made up of one or more pieces joined together
by at least one pixel.

We developed our topology preserving level set method having in view
problems of shape design. There, one specifies in advance certain
conditions about how small, thin, or close certain features of the shape
can get, and then one uses a grid as fine as needed to resolve the details
of the optimal shape. Thus, our method will be different from \cite{han} by
the fact that our method is grid size independent.

\section{Topology-preserving level set method}

A typical shape optimization problem is as follows.  We are given
a cost function $F$ which depends on geometry of the unknown
shape.  The problem is to find a shape such that the cost function
is minimized (at least locally).

Let us represent the shape $D$ as
$$
D=\{x:\phi(x)>0\}.
$$
The optimization problem we wish to solve is
$$
\min_\phi F(\phi),
$$
subject to geometrical and topological constraints on $D$.
The latter constraints are:
\begin{itemize}
\item {\bf Shape topology.} The domain
we design for must have, for example, a fixed number of connected components
or holes.
\item {\bf Component size.} A lower bound
on the size of each component or hole is prescribed.
\item {\bf Distance between components.} A lower bound on the distance between
components or holes is prescribed.  In the case of holes, we
also prescribe a lower bound on the distance from each hole to the 
external boundary of the domain.
\end{itemize}
These  constraints arise naturally in optimal design problems as
we will illustrate in two numerical examples.

It turns out that all these constraints can be handled in a single penalty
formulation.  We will restrict our attention to 2-D problems, even
though the same ideas will work in higher dimensions. 

\begin{figure}
\begin{center}
\includegraphics[width=0.5\textwidth]{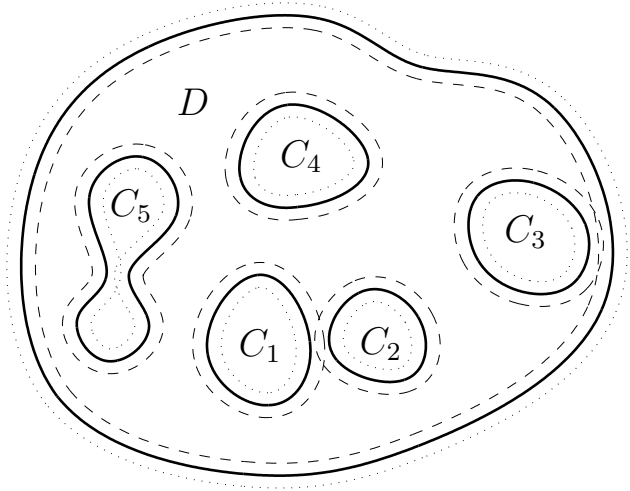} 
\caption{The set $I_d$ (dashed curves) and $E_l$ (dotted curves).}\label{fig1}
\end{center}
\end{figure}

Assume for simplicity that $D$ is a bounded and connected set in $\mathbb R^2$ with
a set of holes inside of it, which are connected components of $\mathbb
R^2\backslash D.$ If $d>0$ and $l>0$ are real numbers, denote 
$$
  I_d=\{x+d\nabla\phi(x): x\in\partial D\},
$$
and
$$
  E_l=\{x-l\nabla\phi(x): x\in\partial D\}.
$$
It follows from \eqref{2} that for $d$ and $l$ small enough, $I_d$ and $E_l$ are made up of points at
distance $d$ and $l$ respectively from $\partial D$.  In fact, for $d=l$,
the union of these two sets is exactly the set of \emph{all} points at
distance $d$ from $\partial D$.  Note that if any two components of
$\mathbb R^2\backslash D$ (we consider the unbounded component too) are
at distance more than $d$ from each other, then $I_d$ is entirely inside
of $D,$ and thus $\phi(x)>0$ on $I_d$.  Also, if the gaps in $D$ are not
too ``small'' or too ``thin'', then $E_l$ is a subset of $\mathbb
R^2\backslash D$, and so $\phi(x)<0$ on $E_l$.

Then, we claim, and using a little bit of geometric intuition it is easy to
see that it is so, that for $d$ and $l$ small numbers, the conditions
$$
  \phi\big(x+d\nabla\phi(x)\big)>0 \mbox{ and }  \phi\big(x-l\nabla\phi(x)\big)<0 \mbox{
    for } x\in\partial D
$$
are a reasonable way of guaranteeing that the holes in $D$ will not
merge, split,  or become too small.  In fact, since the outer boundary of $D$ is
defined by the same level set function, the above also ensures
that the holes will never get too close to the boundary.  These properties
also guarantee that if we start an iterative process with the desired
topology, the iterations cannot change the topology of $D$ as it is
updated.  Thus, these two conditions on the level set function
$\phi$ achieve the constraints of the problem.

To incorporate these conditions into the optimization problem we use
the \emph{logarithmic barrier method}, see
\cite{nocedal-wright-book}. Instead of trying to minimize $F(\phi)$,
consider the problem of minimizing
$F_\varepsilon(\phi)=F(\phi)+\varepsilon H(\phi)$ for $\varepsilon\ll
1,$ where
$$
  H(\phi)=-\int\limits_{\partial D}^{}\!\log \big[\phi\big(x+d\nabla\phi(x)\big)\big]\,ds
   -\int\limits_{\partial D}^{} \log \big[- \phi\big(x-l\nabla\phi(x)\big)\big]\,ds.
$$

To obtain $\phi$ minimizing $F_\varepsilon(\phi)$ we will use the
\emph{steepest descent method}. It amounts to finding the derivative of
$F_\varepsilon(\phi),$ and at each iteration taking a step in the
direction in which the derivative decreases fastest.

In order to calculate the derivative of $F_\varepsilon(\phi)$ we need
the derivatives of $F(\phi)$ and $H(\phi)$. Let $h:\mathbb R^2\to
\mathbb R$ be a test function. For $t$ a real number, $|t|\ll 1,$
$F(\phi+t h)$ will depend on the values of $h$ only close to the boundary of
$D$, as $F$ is a function of the level set of $\phi+t h$, and the way
this level set depends on $h$ was discussed above. We deduce that
$$
  D_\phi F(\phi) \cdot h = \frac{d\,F(\phi+t h)}{dt}\bigg|_{t=0}
$$ 
will only be a function of the restriction of $h$ to $\partial D$. In
many important applications, see \cite{osher-santosa}, it has the form 
\begin{equation}\label{3}
   D_\phi F(\phi) \cdot h = \int\limits_{\partial D}^{}\!U(x)h(x)\,ds,
\end{equation}
for some function $U$ which of course depends on $\phi$ and which can be
calculated numerically. 

The derivative of $H(\phi)$ can be calculated explicitly. Consider a
parameterization $x(s)$ of $\partial D$, with $x'(s)$ having unit norm
for all $s$. $H(\phi+th)$ will be a sum of two integrals over the set
$\{x:(\phi+t h)(x)=0\},$ which, if \eqref{2} holds, is approximately parameterized by $x-t h
(x) \nabla \phi\big(x),$ with $x=x(s)$. One can then find that the derivative
of the first integral in $H(\phi+th)$ at $t=0$ is
\begin{multline}\label{4}
  \int\limits_{\partial D}^{}\!\bigg\{\frac{\big[\nabla\phi(x) h (x)+d  h (x) \nabla^2\phi(x) 
    \nabla\phi(x)-d \nabla h (x)\big]\cdot\nabla\phi\big(x+d\nabla\phi(x)\big)- h\big(x+d\nabla\phi(x)\big)}
{\phi\big(x+d\nabla\phi(x)\big)}\\
+ \log \big[\phi\big(x+d\nabla\phi(x)\big)\big] x' \cdot \big[(\nabla
h(x)\cdot x') \nabla\phi(x) + \big(\nabla ^2 \phi (x)  x'\big)h(x) \big]
    \bigg\} \,ds.
\end{multline}
A similar equality holds for the second term in $H(\phi)$.

Beside the obvious complexity of this expression, note that unlike the
case of $F(\phi),$ this derivative will no longer depend on the values
of the test function $h$ only on $\partial D$. We will make several
approximations. Recall that the purpose of $H(\phi)$ is to
make sure at every step in the optimization process the domain $D$ has
the topology preserved. $H(\phi)$ will grow large only when
$D$ is close to violating the restrictions imposed on it. As far as the
first integral in $H(\phi)$ is concerned, this happens when
$\phi\big(x+d\nabla\phi(x)\big)$ becomes close to zero. Then, the term on
the first line of \eqref{4} is much larger than the second. We will
ignore the term on the second line. Also, on the first line, we have
$\nabla^2\phi(x) \nabla\phi(x)=0,$ which follows from \eqref{2}. In
addition, we will ignore the quantities $-d \nabla h (x)\cdot\nabla\phi\big(x+d\nabla\phi(x)\big)$ and
$-h\big(x+d\nabla\phi(x)\big)$.  We obtain the more manageable expression
$$
\int\limits_{\partial D}^{}\!U_1(x) h(x)\,ds,  
$$
with 
\begin{equation}\label{5}
  U_1(x)=\frac{\nabla\phi(x) \cdot    \nabla\phi\big(x+d\nabla\phi(x)\big)}{\phi\big(x+d\nabla\phi(x)\big)},\, 
  x\in \partial D.
\end{equation}

The derivative of the second integral in $H$ can be
calculated, and then approximated, in the same way. Make the notation
\begin{equation}\label{6}
  U_2(x)=\frac{\nabla\phi(x) \cdot    \nabla\phi\big(x-l\nabla\phi(x)\big)}{\phi\big(x-l\nabla\phi(x)\big)},\, 
  x\in \partial D.
\end{equation}
We obtain 
$$
 D_\phi H(\phi) \cdot h =  \int\limits_{\partial D}^{}\!\big[U_1(x)+U_2(x)\big] h(x)\,ds.
$$

This gives us the following approximate equality 
$$
D_\phi F_\varepsilon(\phi) \cdot h=
  \int\limits_{\partial
   D}^{}\!\big[U(x) + \varepsilon U_1(x)+\varepsilon U_2(x)\big] h(x)\,ds.
$$
If this were an exact equality, the steepest descent
direction for $F_\varepsilon(\phi)$ at $\phi$ would be 
\begin{equation}\label{7}
  u(x)=-[U(x) + \varepsilon U_1(x)+\varepsilon U_2(x)\big],
\end{equation}
where $x\in\partial D$. 
This quantity can be extended continuously to a neighborhood of 
$\partial D$
in the following manner: for $x\in \mathbb R^2$ close to $\partial D$ let
$\widetilde x\in \partial D$ be the unique point such that $\mbox{dist}(x, \partial
D)=\mbox{dist}(x, \widetilde x),$ and set 
\begin{equation}\label{8}
u(x)=\phi(x)+u(\widetilde x).
\end{equation}
Then the
next iteration for $\phi$ would be $\phi+\alpha u,$ where $\alpha>0$ is
the length of the step to be taken in the direction $u.$

But the obtained $u$ is an approximation. It will then clearly not be
the steepest descent direction for $F_\varepsilon(\phi)$. One could
question if it would be a descent direction at all, that is, whether
$F_\varepsilon(\phi)$ would decrease if $\phi$ is replaced by
$\phi+\alpha u.$ After a numerical study we can say that the answer
is no; $F_\varepsilon(\phi)$ could even increase in the
process. Nevertheless, we will argue below that this iterative process
does its job at maintaining the topology constraints. And as far as the
problem of minimizing $F(\phi),$ it is clear that the iterative process
we suggest will give us a sufficiently good approximation to the point
of minimization $\phi$, provided that $\varepsilon$ is small enough.

We will show that, if the level set function $\phi$ is such that two
components of $\{x:\phi(x)<0\}$ are at distance slightly more than $d$
from one another, then $u$ will act as a repelling force, and in
consequence, the components of  $\{x:(\phi+\alpha u)(x)<0\}$ will be
further apart.

\begin{figure}[h]
\begin{center}
\includegraphics[width=0.5\textwidth]{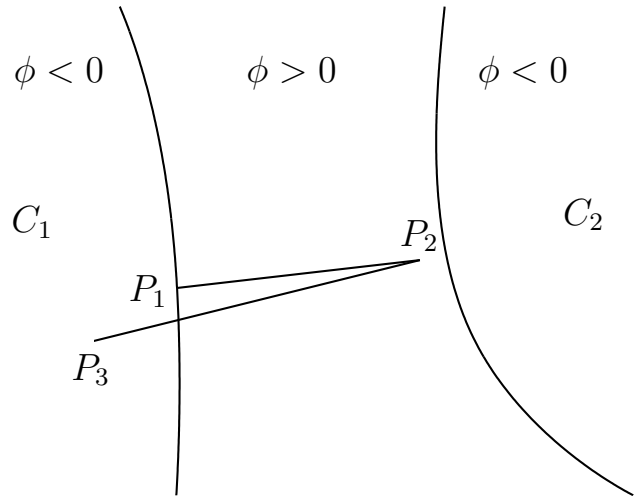} 
\caption{A blow-up of Fig. \ref{fig1}.}\label{fig2}
\end{center}
\end{figure}

Indeed, consider such a situation in Fig. \ref{fig2}. Let $x\in
\partial C_1$ be a point, which in this figure we will denote by $P_1$,
such that $\mbox{dist}(P_1, \partial C_2)$ is slightly larger than $d$.  Let $P_2$ be the point
$x+d\nabla \phi(x)$. Then $P_2$ will be very close to $\partial C_2.$ We
will have $\phi\big(x+d\nabla\phi(x)\big)>0$ but very small. It is easy
to show, and geometrically clear, that the gradient of $\phi$ at $P_2,$
$\nabla\phi\big(x+d\nabla\phi(x)\big)$, which in the figure is
represented by the vector $\stackrel{\longrightarrow}{P_2P_3},$ will
point almost in the opposite direction of
$\stackrel{\longrightarrow}{P_1P_2}$ which is $d\nabla \phi(x)$.  We
find that $U_1(x)$ will be negative and large in magnitude.

Moreover, when the distance between $\partial C_1$ and $\partial C_2$ is
close enough to $d,$ $\varepsilon U_1(x)$ will be larger in
magnitude than $U(x)+\varepsilon U_2(x).$ In consequence, $u(x)$ defined
by \eqref{7} will be positive. Therefore, we have $\phi(x)=0,$ but
$(\phi+\alpha u)(x)>0.$ The same reasoning applies for points $x\in
\partial C_2$ close to $\partial C_1$. This shows that the connected
components of $(\phi+\alpha u)(x)<0$ will be further apart.

It can be argued in the same manner that should a component of $\{x:
\phi(x)<0\}$ get too ``thin'' or too ``small'', then $U_2(x)$ will
serve as a counterweight, forcing it to get ``fatter''. 

Let us note that in order for the above to work, each step size should
not be too big. If the boundary of $D$ moves by more than $\frac{d}{2}$ at some
step, then two components of $\{x: \phi(x)<0\}$ which were at distance
slightly more than $d$ can end up merging without
the penalty functional noticing that. Or, if the boundary moves by more
than $\frac{l}{2},$ a connected component slightly thinner or larger than $l$ might end up splitting or disappearing.
Therefore, at each step one needs to make sure that
\begin{equation}\label{9}
\alpha \max_{x\in \partial D}|u(x)| < K\min(d, l)
\end{equation}
for some $K>0$, as the quantity on the left determines by how much the
boundary of $D$ gets shifted at the given step. Theoretically $K$ can
be allowed to be as large as $\frac{1}{2}$, but since we use a finite grid
size we have to be more conservative. A value of $K=\frac{1}{4}$ works
in practice.

But enforcing \eqref{9} is not enough to guarantee  our
geometrical and topological constraints. The penalty functional
$H(\phi)$ is supposed to take care of this, but it is clear that the
smaller $\varepsilon$ is, the weaker the
influence of $H(\phi)$ in $F_\varepsilon(\phi)$ will be,  and the closer to 
violating the constraints $\phi$ will get, before this penalty functional kicks in. Thus, at each
iteration one needs to first take a step size $\alpha$ satisfying
\eqref{9}, and still check after updating $\phi$ to $\phi+\alpha
u$ whether $H(\phi)$ is
defined. If not, one needs to decrease the step
size $\alpha$, for example by halving it, until $H(\phi)$ is defined. If
no amount of decreasing $\alpha$ helps, one needs to either increase
$\varepsilon$ or decrease the
grid size, and restart the algorithm.

A pseudo-code for the algorithm is as follows.

\begin{center}
\fbox{
\begin{minipage}{6in}
{\tt
initial guess for $\phi$\\
do while not optimal\\
\mbox{ }$\bullet$ compute the descent direction $u$ (use \eqref{3},
\eqref{5}, \eqref{6}, \eqref{7}, and \eqref{8})\\
\mbox{ }$\bullet$ choose a step size $\alpha$ satisfying \eqref{9} for
which $H(\phi+\alpha u)$ is defined\\
\mbox{ }$\bullet$ update $\phi$ to $\phi+\alpha u$\\
\mbox{ }$\bullet$ reinitialize $\phi$ to satisfy \eqref{1}
}
\end{minipage}
}
\end{center}

We note that if at some point the contour $\{x:\phi(x)=0\}$ develops
sharp angles, then the functional $H(\phi)$ might not be defined (this can be
seen from Fig. \ref{fig1}). To prevent this from happening, one can
smooth $\phi$ a bit at each iteration. For $\phi$ discretized on a square grid we used the
procedure 
$$
\phi_{i,j}\to\frac{\phi_{i,j}+\phi_{i-1,j}+\phi_{i+1,j}+\phi_{i,j-1}+\phi_{i,j+1}}{5}.
$$

Also, for fine grids it becomes expensive to reinitialize $\phi$ according to \eqref{1}.
To make this computation faster we reinitialized $\phi$ only in a neighborhood of the set $\{x:\phi(x)=0\}.$
For more performance one could use the fast re-distancing algorithms suggested in 
\cite{smereka, strain, sussman}.

Lastly, sometimes one might wish to introduce additional constraints of
the form $G(\phi)=\mbox{const.}$ in the optimization problem. An example of such a constraint is the
requirement that the area of the set $\{x:\phi(x)>0\}$ be kept fixed, which
we will use in the two numerical examples below. Then one needs to
modify the descent direction $u$ as described in \cite{osher-santosa}.

\section{Numerical examples}

In the first example, we consider the problem of finding a domain that has
the smallest perimeter, subject to the constraint that the area of the domain being fixed.
Thus, the functional to minimize is 
$$
  F(\phi)=\int\limits_{\{\phi=0\}}^{}\!1\,ds,
$$
with the constraint
$$
  G(\phi)=\int\limits_{\{\phi>0\}}^{}\!1\,dx\,dy=\mbox{const.}
$$
The starting shape is a region with seven subdomains, each one an
ellipse with aspect ratio 1.3, as shown in Fig. \ref{fig3} on the left.  The
center ellipse has a slightly bigger (20\%) size than the rest.  The
distance between the centers of the ellipses is 4, and the smallest
semi-axis of the surrounding ellipses is 1.

\def\width{0.4}
\begin{figure}[t]
\begin{center}
\fbox{\includegraphics[clip, width=\width\textwidth]{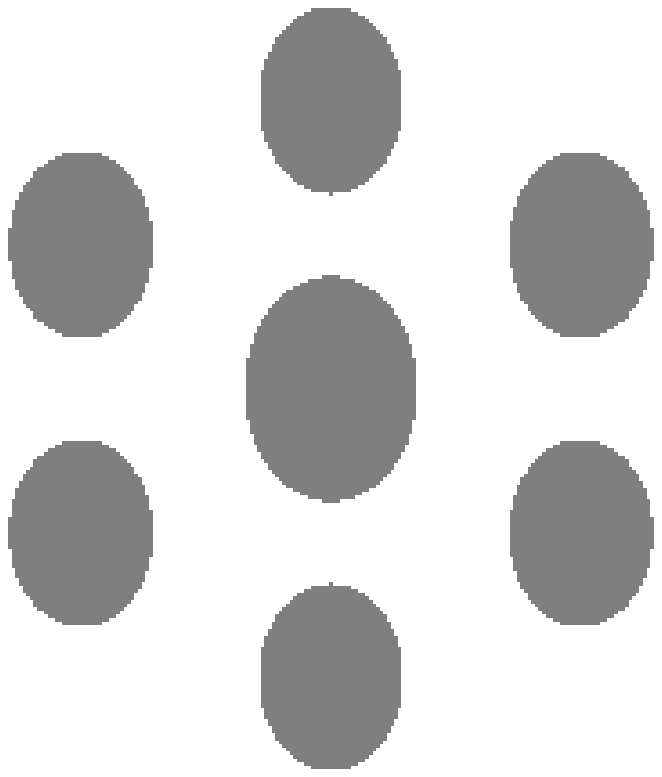}} \hspace{0.5cm}
\fbox{\includegraphics[clip, width=\width\textwidth]{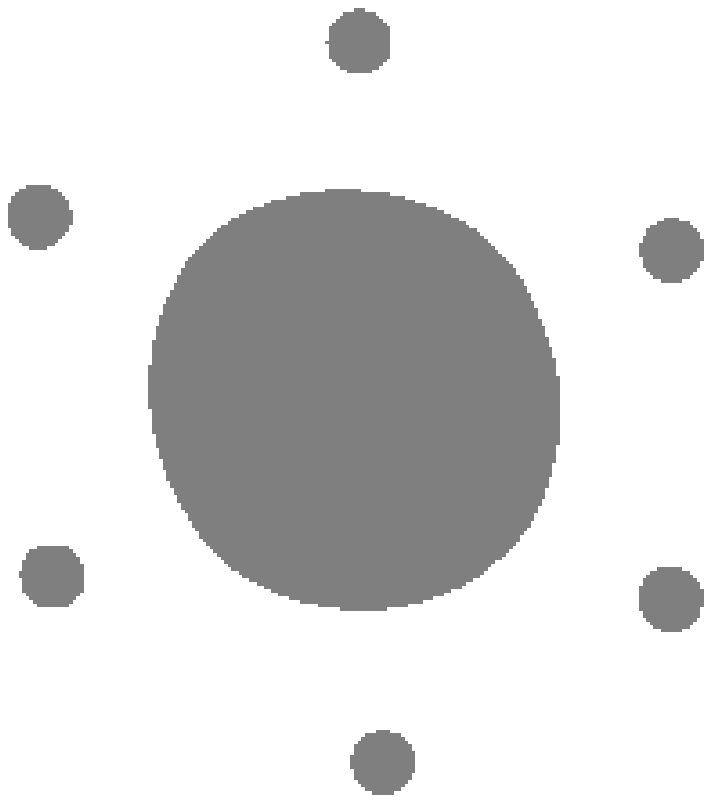}}
\caption{The initial and optimized shape for example 1.}\label{fig3}
\end{center}
\end{figure}

If we do not constrain the
topology or geometry, the optimal solution would be a circle whose area
is equal to the area of the original seven subdomains.
If we do enforce these constraints, minimizing instead the functional 
$$
F(\phi) + \varepsilon H(\phi),
$$
we obtain the picture in Fig. \ref{fig3} on the right.

For this calculation we set $d=l=0.8$, $\varepsilon=0.2$ and consider a square grid of
size $h=0.05$ (each square is further split into two triangles, to make
it easier to keep track of the set $\{x:\phi(x)=0\}$). 

We find that the ``satellite'' components of the central domain do not
disappear, but became of size slightly larger than $l.$ 

We note that that
the resulting large domain in the center is not perfectly circular. This
because the steepest descent direction for $F(\phi)$ will be $\Delta
\phi=\phi_{xx}+\phi_{yy}.$ We need to calculate this quantity
numerically, and after reinitializing $\phi$ according to \eqref{1} it
is not smooth enough for $\Delta\phi$ to be calculated
accurately. Smoothing $\phi$ as noted in the previous section helped a
bit, this is how this picture was obtained.  We found that if we perform
additional smoothing then the result in Fig. \ref{fig2} will look
more circular.  This artifact does not show up in the next example, as
then one does not need to calculate second-order derivatives of $\phi.$

\def\width{0.4}
\begin{figure}[t]
\begin{center}
\fbox{\includegraphics[clip, width=\width\textwidth]{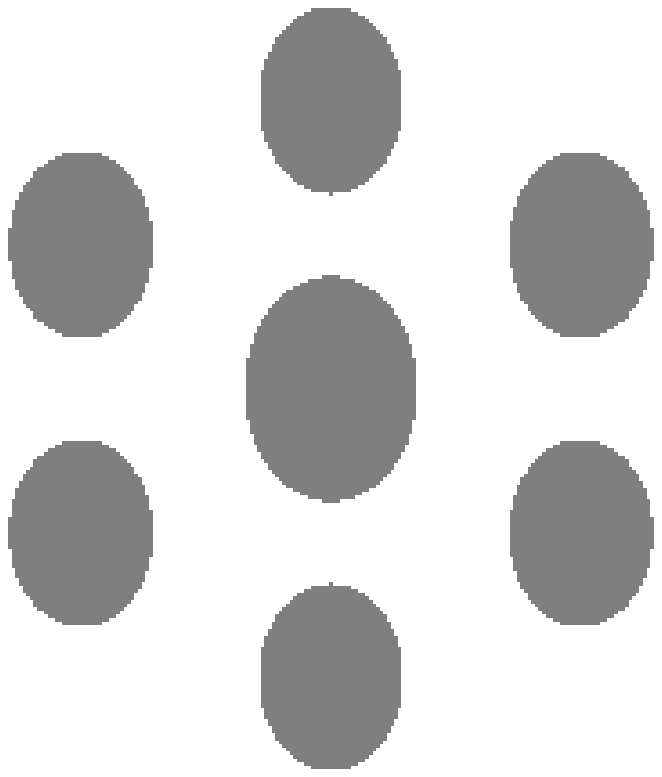}} \hspace{0.5cm}
\fbox{\includegraphics[clip, width=\width\textwidth]{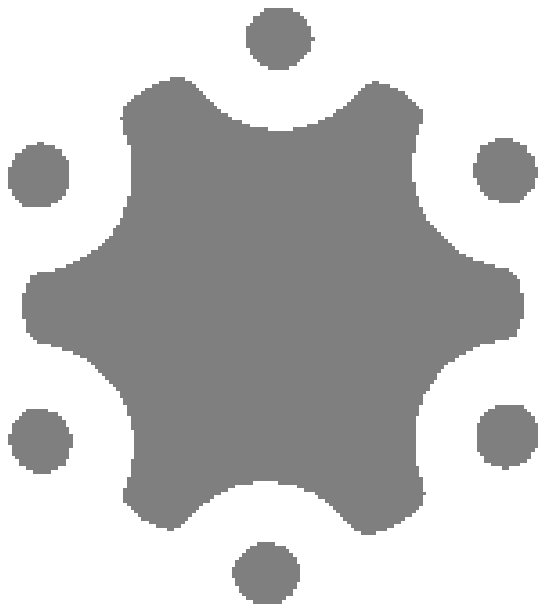}}
\caption{The initial and optimized shape for example 2.}\label{fig4}
\end{center}
\end{figure}

In the second example we examine the problem of 
minimizing the functional
$$
  F(\phi)=\int\limits_{\{\phi>0\}}^{}\!(x^2+y^2)\,dx\,dy.
$$
We again enforce the area constraint $G(\phi)=\mbox{const.},$ and we use the same
values for $d,$ $l$ and $h$. We set $\varepsilon=0.4$.  (The value of
$\varepsilon$ which is relatively small, and in the same time be not
small enough that the algorithm fails converge for a given grid size is
determined by trial and error, and it depends on the problem.) In
absence of topological constraints, these seven ellipses would merge to
form a large circle.  The topological constraints prevent them from
doing so, as we see from Fig. \ref{fig4}.

\section{Discussion}
In this paper we introduced a penalty functional which makes it
possible to use the level set method in problems with topology and geometry
constraints.  Our method allows for topological constraints independent
of the grid size (that is, for given $d$ and $l$, the grid size $h$ can
be chosen as small as desired), which is a key difference with the
method suggested in \cite{han}.

\bibliographystyle{plain}
\bibliography{logbar.bib}
\end{document}